\documentclass[10pt, reqno]{amsart}

\usepackage{amssymb}
\usepackage{enumitem}
\usepackage{hyperref}
\usepackage[noabbrev]{cleveref}
\usepackage{mathtools}
\usepackage{microtype}[final]
\usepackage{quiver}
\usepackage{subcaption}
\usepackage[format=plain,
            labelfont={bf},
            textfont=it]{caption}
\usepackage{todonotes}

\usepackage[backend=biber, doi=false, giveninits=true, isbn=false, maxbibnames=99,
sortcites=true, sorting=nyt, style=numeric-comp, url=false]{biblatex} 
\renewbibmacro{in:}{}  
\allowdisplaybreaks

\numberwithin{equation}{section}

    \renewcommand{\epsilon}{\varepsilon}
    \renewcommand{\phi}{\varphi}

    \newcommand{\defeq}{\vcentcolon=}
    \newcommand{\eqdef}{=\vcentcolon}

    \newcommand{\set}[1]{\{ #1 \}}

    \newcommand{\R}{\mathbb{R}}

    \newcommand{\abs}[1]{\lvert #1 \rvert}

    \newcommand{\inner}[2]{\langle #1, #2 \rangle}

    \newcommand{\norm}[1]{\lVert #1 \rVert}
    \newcommand{\Norm}[1]{\left\lVert #1 \right\rVert}

    \newcommand{\les}{\lesssim}

    \newcommand{\mes}{\abs}

    \newcommand{\BL}{\mathrm{BL}}
    \newcommand{\BLG}{\mathrm{BL}_\mathrm{G}}
    \newcommand{\BLCD}{\mathrm{BLCD}}
    \newcommand{\BLCDG}{\mathrm{BLCD}_\mathrm{G}}

    \newcommand{\A}{\mathcal{A}}

    \renewcommand{\vec}{\mathbf}
    
    \newcommand{\subspeq}{\leq}
    \newcommand{\tp}{\top}
    \newcommand{\card}[1]{\#(#1)}
    \newcommand{\codim}{\operatorname{codim}}
    \newcommand{\spn}{\operatorname{span}}

\newtheorem{theorem}{Theorem}[section]

\newtheorem{corollary}[theorem]{Corollary}

\theoremstyle{definition}
\newtheorem{definition}[theorem]{Definition}
\newtheorem{example}[theorem]{Example}
\newtheorem*{ack}{Acknowledgments}

\theoremstyle{remark}
\newtheorem{remark}[theorem]{Remark}

\begin{document}

\title{Quiver Brascamp--Lieb Inequalities}

\author[N.~Hu]{Nicholas Hu}
\address{Department of Mathematics, UCLA}
\email{njhu@math.ucla.edu}

\begin{abstract} 
    We formulate generalized Brascamp--Lieb inequalities for representations of bipartite
    quivers and establish necessary and sufficient conditions for such inequalities.
    Notably, we show contra Lieb that Gaussians do not saturate certain types of quiver
    Brascamp--Lieb inequalities.
\end{abstract}

\maketitle

\section{Introduction}

The Brascamp--Lieb inequalities \cite{BrascampLieb} are an important family of inequalities
in analysis that subsumes several inequalities significant in their own right, including
H\"older's inequality, the Loomis--Whitney inequality, and Young's inequality. Several
variants and extensions of these inequalities have been developed, such as reverse
Brascamp--Lieb inequalities \cite{Barthe}, perturbed Brascamp--Lieb inequalities
\cite{Bennett3}, and adjoint Brascamp--Lieb inequalities \cite{Bennett4}, some of which have
proved to be very useful in Fourier restriction theory \cite{Zhang}. 

In addition, algorithms have been devised to compute optimal constants in Brascamp--Lieb
inequalities \cite{Garg,Weber}. The algorithm of Garg, Gurvits, Oliveira, and Wigderson
\cite{Garg} does so by relating these constants to the so-called \emph{capacity} of
completely positive operators, and this notion of capacity was generalized by Chindris and
Derksen \cite{Chindris} to algebraic objects known as \emph{quiver representations}. In this
paper, we come full circle by formulating and studying Brascamp--Lieb inequalities for
quivers.

\subsection{Brascamp--Lieb inequalities}

Let us begin by briefly reviewing the theory of ordinary Brascamp--Lieb inequalities.

\begin{definition}[Brascamp--Lieb inequality]
    Let $H$ and $H^1, \dots, H^m$ be nontrivial finite-dimensional Hilbert spaces, and for
    each $1 \leq j \leq m$, let $B_j : H \to H^j$ be a surjective linear map and $p_j \in
    [1, \infty]$. A \emph{Brascamp--Lieb inequality} is an inequality of the form 
    \begin{equation} \label{bl}
        \int_H \prod_{j=1}^m f_j \circ B_j \, dx \les \prod_{j=1}^m \norm{f_j}_{L^{p_j}(H^j)}
    \end{equation} 
    that holds for all measurable functions $f_j : H^j \to [0, \infty]$, where the implicit
    constant (see \Cref{S:notation}) is allowed to be infinite.
\end{definition}

Of course, a Brascamp--Lieb inequality is only useful when the implicit constant is finite.
The conditions under which this occurs are well known \cite{Bennett,Bennett2}.

\begin{theorem} \label{T:bl}
    Let $\BL(\vec{B}, \vec{p})$ denote the smallest constant for which inequality
    \labelcref{bl} holds, where $(\vec{B}, \vec{p}) \defeq ((B_1, \dots, B_m), \allowbreak
    (p_1, \dots, p_m))$; it is known as the \emph{Brascamp--Lieb constant} for the
    \emph{Brascamp--Lieb datum} $(\vec{B}, \vec{p})$. Then $\BL(\vec{B}, \vec{p})$ is finite
    if and only if the \emph{scaling condition}
    \begin{equation} \label{scal}
        \dim(H) = \sum_{j=1}^m p_j^{-1} \dim(H^j)
    \end{equation}
    and the \emph{dimension condition}
    \begin{equation} \label{dim}
        \dim(V) \leq \sum_{j=1}^m p_j^{-1} \dim(B_j V) \quad \text{for all $V \subspeq H$}
    \end{equation}
    hold. \textup(Here $V \subspeq H$ means that $V$ is a subspace of $H$.\textup)
\end{theorem}

Furthermore, Lieb \cite{Lieb} showed that the optimal constant is unchanged when the
functions are restricted to be Gaussians of the form $f_j(x) \defeq e^{-\pi \inner{A_j x}{x}
p_j^{-1}}$ for each $1 \leq j \leq m$, 
where $A_j \succ 0$ with respect to the Loewner order. Since
\begin{equation*}
    \int_{H^j} e^{-\pi \inner{A_j x}{x}} \, dx = \det(A_j)^{-1/2},
\end{equation*}
this yields the following result and formula for Brascamp--Lieb constants.

\begin{theorem} \label{T:lieb}
    Let $\BLG(\vec{B}, \vec{p})$ denote the smallest constant for which inequality
    \labelcref{bl} holds when $f_j(x) \defeq e^{-\pi \inner{A_j x}{x} p_j^{-1}}$ for some $A_j
    \succ 0$ \textup(for each $1 \leq j \leq m$\textup). Then
    \begin{equation} \label{lieb}
        \BL(\vec{B}, \vec{p}) 
        = \BLG(\vec{B}, \vec{p})
        = \sup_{A_j \succ 0} \left[\frac{\prod_{j=1}^m
        \det(A_j)^{p_j^{-1}}}{\det\left(\sum_{j=1}^m p_j^{-1} B_j^\tp A_j
        B_j\right)}\right]^{1/2}.
    \end{equation}
    \textup(Here the supremum is taken over $A_1 \succ 0, \dots, A_m \succ 0$.\textup)
\end{theorem}

\begin{example}[Young's convolution inequality]
    Let $H \defeq \R^d \times \R^d$ and $H^j \defeq \R^d$ for $1 \leq j \leq 3$ 
    (with their standard inner products), 
    and let $B_1(x, y)
    \defeq x$, $B_2(x, y) \defeq y$, and $B_3(x, y) \defeq x-y$ for all $(x, y) \in H$. 
    Then inequality \labelcref{bl} is \emph{Young's convolution inequality}
    \[
        \int_{\R^d} \int_{\R^d} f_1(x) f_2(y) f_3(x-y) \, dx \, dy
        \les
        \norm{f_1}_{L^{p_1}(\R^d)}
        \norm{f_2}_{L^{p_2}(\R^d)}
        \norm{f_3}_{L^{p_3}(\R^d)}.
    \]
    It is a result of Beckner \cite{Beckner} and Brascamp and Lieb \cite{BrascampLieb} that
    \[
        \BL(\vec{B}, \vec{p}) = \left(\prod_{j=1}^3 \frac{p_j^{1/p_j}}{{p_j'}^{1/p_j'}}\right)^{d/2}
        \quad
        \text{if $\sum_{j=1}^3 \frac{1}{p_j} = 2$},
    \]
    where $1/p_j + 1/p_j' = 1$,
    and $\BL(\vec{B}, \vec{p}) = \infty$ for all other $\vec{p}$ (under the usual
    interpretations of $1/\infty = 0$ and $\infty^{1/\infty} = 1$).
\end{example}

\subsection{Quiver Brascamp--Lieb inequalities}

The relationships between the spaces and the linear maps in inequality \labelcref{bl} 
can be depicted by a type of graph known as a \emph{quiver}.
\begin{definition}[Quiver]
    A \emph{quiver} (or \emph{directed multigraph}) $\mathcal{Q}$ is a tuple $(\mathcal{V},
    \A, s, t)$ consisting of a set $\mathcal{V}$ of \emph{vertices}, a set
    $\A$ of \emph{arrows} (or \emph{edges}), a function $s : \A \to
    \mathcal{V}$ specifying the \emph{source} (or \emph{tail}) of each arrow, and a function
    $t : \A \to \mathcal{V}$ specifying the \emph{target} (or \emph{head}) of each
    arrow. If $\mathcal{V}$ is the disjoint union of two sets $\mathcal{V}_1$ and
    $\mathcal{V}_2$ with $s(\A) \subseteq \mathcal{V}_1$ and $t(\A)
    \subseteq \mathcal{V}_2$, the quiver is said to be \emph{bipartite}.
\end{definition}
\begin{definition}[Quiver representation]
    A \emph{representation} of a quiver is an assignment of a vector space to each vertex
    and a linear map to each arrow.
\end{definition}
In what follows, we will only consider quivers with finitely many vertices and arrows and
will generally identify each quiver with a specific representation of it by nontrivial
finite-dimensional Hilbert spaces and surjective linear maps (so, strictly speaking, the
sets of vertices and arrows will be multisets). For instance, the quiver corresponding to
inequality \labelcref{bl} is the bipartite \emph{$m$-subspace quiver} defined by
$\mathcal{V} = \set{H} \cup \set{H^j}_{j=1}^m$, $\A = \set{B_1, \dots, B_m}$,
$s(B_j) = H$, and $t(B_j) = H^j$, as illustrated in \Cref{F:subsp-quiv}.

\begin{figure}[htbp]
    \centering
    \begin{subfigure}{0.5\textwidth}
        \centering
        \begin{tikzcd}[row sep=tiny]
            && {H^1} \\
            && {H^2} \\
            H && \vdots \\
            && \vdots \\
            && {H^m}
            \arrow["{B_1}", from=3-1, to=1-3]
            \arrow["{B_2}"', from=3-1, to=2-3]
            \arrow["{B_m}"', from=3-1, to=5-3]
        \end{tikzcd}
        \caption{The $m$-subspace quiver}
        \label{F:subsp-quiv}
    \end{subfigure}%
    \begin{subfigure}{0.5\textwidth}
        \centering
        \begin{tikzcd}[row sep=small]
            {H_1} && {H^1} \\
            {H_2} && {H^2} \\
            \vdots && \vdots \\
            {H_n} && {H^m} \\
            && {}
            \arrow[curve={height=-12pt}, from=1-1, to=1-3]
            \arrow[from=1-1, to=1-3]
            \arrow[from=2-1, to=1-3]
            \arrow[from=2-1, to=4-3]
            \arrow[from=4-1, to=2-3]
        \end{tikzcd}
        \caption{A bipartite quiver}
        \label{F:bipart-quiv}
    \end{subfigure}

    \caption{}
\end{figure}

Recently, Chindris and Derksen \cite{Chindris} considered a quantity analogous to that in
formula \labelcref{lieb} for a general bipartite quiver with $\mathcal{V} =
\set{H_i}_{i=1}^n \cup \set{H^j}_{j=1}^m$, $s(\A) = \set{H_i}_{i=1}^n$, and
$t(\A) = \set{H^j}_{j=1}^m$; see \Cref{F:bipart-quiv}.  Namely, if
$\A_{ij}$ denotes the set of arrows from $H_i$ to $H^j$, $B_a$ denotes the
surjective linear map representing an arrow $a \in \A$, and $p_j \in [1, \infty]$
for each $1 \leq j \leq m$, it was shown that the quantity
\begin{equation} \label{qsup}
    \sup_{A_j \succ 0} \left[\frac{\prod_{j=1}^m
    \det(A_j)^{p_j^{-1}}}{\prod_{i=1}^n \det\left(\sum_{j=1}^m \sum_{a \in \A_{ij}} p_j^{-1} B_a^\tp A_j
    B_a\right)}\right]^{1/2}
\end{equation}
is finite if and only if the scaling condition
\begin{equation} \label{qscal}
    \sum_{i=1}^n \dim(H_i) = \sum_{j=1}^m p_j^{-1} \dim(H^j)
\end{equation}
and the dimension condition
\begin{equation} \label{qdim}
    \sum_{i=1}^n \dim(V_i) \leq 
    \sum_{j=1}^m p_j^{-1} \dim\left(\sum_{i=1}^n \sum_{a \in \A_{ij}} B_a V_i\right) 
    \quad \text{for all $V_i \subspeq H_i$}
\end{equation}
hold. Clearly, conditions \labelcref{qscal,qdim} and the supremum \labelcref{qsup} reduce to
conditions \labelcref{scal,dim} and the supremum \labelcref{lieb} in the case of the
$m$-subspace quiver.

However, the natural question of whether the supremum \labelcref{qsup} is the optimal
constant in some generalization of the Brascamp--Lieb inequality was left uninvestigated.
Indeed, a ``quiver Brascamp--Lieb inequality'' has yet to be formulated, so our first task
is to reinterpret this algebraic result as an analytic one.
Na\"ively taking $f_j(x) \defeq e^{-\pi \inner{A_j x}{x} p_j^{-1}}$ as before allows one to write
the supremum \labelcref{qsup} as
\begin{equation*}
    \sup_{A_j \succ 0}
    \frac{\prod_{i=1}^n \int_{H_i} \prod_{j=1}^m \prod_{a \in \A_{ij}} f_j \circ
    B_a \, dx}{\prod_{j=1}^m \norm{f_j}_{L^{p_j}(H^j)}},
\end{equation*}
but this misleadingly suggests the untenable inequality
\begin{equation*}
    \prod_{i=1}^n \int_{H_i} \prod_{j=1}^m \prod_{a \in \A_{ij}} f_j \circ
    B_a \, dx \les \prod_{j=1}^m \norm{f_j}_{L^{p_j}(H^j)},
\end{equation*}
which in general has a different number of $f_j$ on each side.
Aiming for the more viable inequality
\begin{equation} \label{blq}
    \prod_{i=1}^n \int_{H_i} \prod_{j=1}^m \prod_{a \in \A_{ij}} f_j \circ
    B_a \, dx \les \prod_{i=1}^n \prod_{j=1}^m \prod_{a \in \A_{ij}} 
    \norm{f_j}_{L^{p_j}(H^j)} = \prod_{j=1}^m \norm{f_j}_{L^{p_j}(H^j)}^{\alpha_j},
\end{equation}
where $\alpha_j \defeq \sum_{i=1}^n \card{\A_{ij}}$ (and $\card{\A_{ij}}$ denotes the
cardinality of $\A_{ij}$), one can verify that scaling $p_j^{-1}$ by $\alpha_j$ in 
the supremum \labelcref{qsup} yields
\begin{equation*}
    C_{\mathcal{Q}, \vec{p}}^{-1/2} \cdot
    \sup_{A_j \succ 0}
    \frac{\prod_{i=1}^n \int_{H_i} \prod_{j=1}^m \prod_{a \in \A_{ij}} f_j \circ
    B_a \, dx}{\prod_{j=1}^m \norm{f_j}_{L^{p_j}(H^j)}^{\alpha_j}}
\end{equation*}
when $f_j(x) \defeq e^{-\pi \inner{A_j x}{x} \alpha_j p_j^{-1}}$, where
\begin{equation*}
    C_{\mathcal{Q}, \vec{p}} \defeq \prod_{j=1}^m \alpha_j^{\alpha_j p_j^{-1} \dim(H^j)}.
\end{equation*}
Thus, the result of Chindris and Derksen amounts to the theorem below.

\begin{theorem} \label{T:blq}
    Let $\BLCDG(\mathcal{Q}, \vec{p})$ denote the optimal constant in inequality
    \labelcref{blq} when $f_j(x) \defeq e^{-\pi \inner{A_j x}{x} \alpha_j p_j^{-1}}$ for some
    $A_j \succ 0$, so that
    \begin{equation*}
        \BLCDG(\mathcal{Q}, \vec{p}) = C_{\mathcal{Q}, \vec{p}}^{1/2} \cdot
        \sup_{A_j \succ 0} \left[\frac{\prod_{j=1}^m
        \det(A_j)^{\alpha_j p_j^{-1}}}{\prod_{i=1}^n \det\left(\sum_{j=1}^m \sum_{a \in
            \A_{ij}} \alpha_j p_j^{-1} B_a^\tp A_j
        B_a\right)}\right]^{1/2}.
    \end{equation*}
    Then $\BLCDG(\mathcal{Q}, \vec{p})$ is finite if and only if the \emph{scaling
    condition}
    \begin{equation} \label{blgqscal}
        \sum_{i=1}^n \dim(H_i) = \sum_{j=1}^m \alpha_j p_j^{-1} \dim(H^j)
    \end{equation}
    and the \emph{dimension condition}
    \begin{equation} \label{blgqdim}
        \sum_{i=1}^n \dim(V_i) \leq 
        \sum_{j=1}^m \alpha_j p_j^{-1} \dim\left(\sum_{i=1}^n \sum_{a \in \A_{ij}} B_a V_i\right) 
        \quad \text{for all $V_i \subspeq H_i$}
    \end{equation}
    hold. 
\end{theorem}

One might hope that a generalization of \Cref{T:lieb} holds so that the optimal constant
$\BLCD(\mathcal{Q}, \vec{p})$ in inequality \labelcref{blq} (for general functions)
coincides with $\BLCDG(\mathcal{Q}, \vec{p})$, or at least that conditions
\labelcref{blgqscal,blgqdim} are also sufficient for the finiteness of $\BLCD(\mathcal{Q},
\vec{p})$. Unfortunately, this turns out to be false -- the finiteness of the general
constant is in fact equivalent to the scaling condition \labelcref{blgqscal} along with the
stronger dimension condition \labelcref{blqadim} given below.  We will prove this by
considering ostensibly more general inequalities in which a function is associated with each
\emph{arrow} instead of each target space.

\begin{definition}[Quiver Brascamp--Lieb inequalities]
    Let $\mathcal{Q}$ be a bipartite quiver represented by nontrivial finite-dim\-ensional
    Hilbert spaces $H_1, \dots, H_n$ and $H^1, \dots,\allowbreak H^m$ and surjective linear
    maps from the former to the latter. In addition, let $\A_{ij}$ denote the set
    of arrows from $H_i$ to $H^j$, $B_a$ denote the map representing an arrow $a$, and $p_j
    \in [1, \infty]$ for each $1 \leq j \leq m$. A \emph{quiver Brascamp--Lieb inequality}
    is an inequality of the form
    \begin{equation} \label{blqa}
        \prod_{i=1}^n \int_{H_i} \prod_{j=1}^m \prod_{a \in \A_{ij}} f_a \circ
        B_a \, dx \les \prod_{i=1}^n \prod_{j=1}^m \prod_{a \in \A_{ij}} 
        \norm{f_a}_{L^{p_j}(H^j)}
    \end{equation} 
    that holds for all measurable functions $f_a : H^j \to [0, \infty]$.
\end{definition}

\begin{example} \label{E:example}
    Consider the quiver $\mathcal{Q}$ in \Cref{F:example}, where $H_1 \defeq \R^3$, $H^1
    \defeq \R^2$, $B_1(x_1, x_2, x_3) \defeq (x_1, x_2)$, and $B_2(x_1, x_2, x_3) \defeq
    (x_2, x_3)$; and let $\vec{p} = (p_1) \defeq (\frac{4}{3})$.
    \begin{figure}[htbp]
        \centering
        \begin{tikzcd}
            {H_1} && {H^1}
            \arrow["{B_1}", curve={height=-6pt}, from=1-1, to=1-3]
            \arrow["{B_2}"', curve={height=6pt}, from=1-1, to=1-3]
        \end{tikzcd}
        \caption{}
        \label{F:example}
    \end{figure}

    Then inequality \labelcref{blq} reads
    \[
        \int_{\R^3} f_1(x_1, x_2) f_1(x_2, x_3) \, dx
        \les \norm{f_1}_{L^{4/3}(\R^2)} \norm{f_1}_{L^{4/3}(\R^2)}
        = \norm{f_1}_{L^{4/3}(\R^2)}^2,
    \]
    whereas inequality \labelcref{blqa} reads
    \[
        \int_{\R^3} f_1(x_1, x_2) f_2(x_2, x_3) \, dx
        \les \norm{f_1}_{L^{4/3}(\R^2)} \norm{f_2}_{L^{4/3}(\R^2)}.
    \]
\end{example}

For quiver Brascamp--Lieb inequalities in the form of inequality \labelcref{blqa}, a
generalization of \Cref{T:bl} holds.

\begin{theorem} \label{T:blqa}
    Let $\BL(\mathcal{Q}, \vec{p})$ denote the optimal constant in inequality
    \labelcref{blqa}. Then $\BL(\mathcal{Q}, \vec{p})$ is finite if and only if the
    \emph{scaling condition}
    \begin{equation} \label{blqascal}
        \sum_{i=1}^n \dim(H_i) = \sum_{i=1}^n \sum_{j=1}^m \sum_{a \in \A_{ij}}
        p_j^{-1} \dim(H^j)
    \end{equation}
    and the \emph{dimension condition}
    \begin{equation} \label{blqadim}
        \sum_{i=1}^n \dim(V_i) \leq 
        \sum_{i=1}^n \sum_{j=1}^m \sum_{a \in \A_{ij}} 
        p_j^{-1} \dim(B_a V_i) 
        \quad \text{for all $V_i \subspeq H_i$}
    \end{equation}
    hold.
\end{theorem}

In actuality, such inequalities are no more general than Chindris--Derksen-type
inequalities, as the following result shows.

\begin{theorem} \label{T:equiv}
    Let $\BLCD(\mathcal{Q}, \vec{p})$ denote the optimal constant in inequality
    \labelcref{blq}. Then
    \[
        \BLCD(\mathcal{Q}, \vec{p}) 
        \leq \BL(\mathcal{Q}, \vec{p}) 
        \leq \prod_{j=1}^m \alpha_j^{\alpha_j} \cdot \BLCD(\mathcal{Q}, \vec{p}).
    \]
    \textup(As above, $\alpha_j = \sum_{i=1}^n \card{\A_{ij}}$.\textup)
\end{theorem}

\begin{corollary} \label{C:equiv}
    Let $\BLCD(\mathcal{Q}, \vec{p})$ denote the optimal constant in inequality
    \labelcref{blq}. Then $\BLCD(\mathcal{Q}, \vec{p})$ is finite if and only if the
    \emph{scaling condition} \labelcref{blqascal}
    \begin{equation*}
        \sum_{i=1}^n \dim(H_i) = \sum_{j=1}^m \alpha_j p_j^{-1} \dim(H^j)
    \end{equation*}
    and the \emph{dimension condition} \labelcref{blqadim}
    \begin{equation*}
        \sum_{i=1}^n \dim(V_i) \leq 
        \sum_{i=1}^n \sum_{j=1}^m \sum_{a \in \A_{ij}} p_j^{-1} \dim(B_a V_i) 
        \quad \text{for all $V_i \subspeq H_i$}
    \end{equation*}
    hold.
\end{corollary}

As a result, we find that Gaussians do not saturate Chindris--Derksen-type inequalities in
general.

\begin{corollary} \label{C:nolieb}
    There exists a quiver Brascamp--Lieb datum $(\mathcal{Q}, \vec{p})$ for which
    $\BLCDG(\mathcal{Q}, \vec{p}) < \BLCD(\mathcal{Q}, \vec{p}) = \infty$.
\end{corollary}

The sufficiency and necessity of the conditions in \Cref{T:blqa} will be separately
established in \Cref{S:conditions}. The proof of sufficiency involves splitting the
bipartite quiver in question into multiple subspace quivers to which \Cref{T:bl} applies;
inequality \labelcref{blqa} is then obtained as the product of inequalities of the form
\labelcref{bl}. The proof of necessity is an adaptation of the scaling argument used by
Bennett, Carbery, Christ, and Tao \cite{Bennett}.

Following this, in \Cref{S:examples}, we will show that every instance of inequality
\labelcref{blq} can be realized as an instance of inequality \labelcref{blqa} and
vice-versa, which will prove \Cref{T:equiv}. \Cref{C:equiv} is merely a restatement of
\Cref{T:blqa} for comparison with \Cref{T:blq}. Finally, we will verify that for
the datum in \Cref{E:example}, the conditions in \Cref{T:blq} hold but those in
\Cref{T:blqa} or \Cref{C:equiv} do not, which entails \Cref{C:nolieb} -- in essence, the
failure of \Cref{T:lieb} for quivers.

\begin{ack}
    We thank Terence Tao for his guidance and suggestions in the research for and writing of
    this paper, as well as Calin Chindriș for his thought-provoking discussions. We are
    obliged to the reviewer for helpful comments.
\end{ack}

\subsection{Notation} \label{S:notation}

We employ the standard notation $A \les B$ to indicate that $A \leq CB$ for some constant $C
> 0$; if $A \les B$ and $B \les A$, we write $A \approx B$.  Occasionally, we adjoin
subscripts to this notation to indicate dependence of the constant $C$ on other parameters;
for instance, we write $A \les_{\alpha, \beta} B$ when $A \leq CB$ for some constant $C > 0$
depending on $\alpha, \beta$.

We also use $\card{{}\cdot{}}$ for the cardinality of a finite set, $\mes{{}\cdot{}}$ for
the measure of a set, and $\norm{{}\cdot{}}$ for the norm of a vector.

\section{Conditions for quiver Brascamp--Lieb inequalities} \label{S:conditions}

First, we prove that the scaling and dimension conditions in \Cref{T:blqa} are sufficient
and necessary.

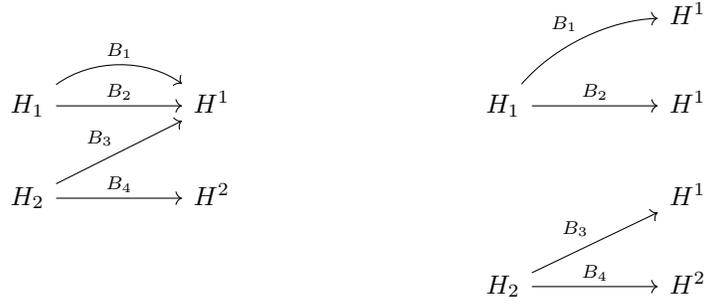
\begin{figure}[htbp]
    \centering
    \begin{subfigure}{0.5\textwidth}
        \centering
        \begin{tikzcd}
            {H_1} && {H^1} \\
            {H_2} && {H^2}
            \arrow["{B_2}", from=1-1, to=1-3]
            \arrow["{B_3}", from=2-1, to=1-3]
            \arrow["{B_4}", from=2-1, to=2-3]
            \arrow["{B_1}", curve={height=-18pt}, from=1-1, to=1-3]
        \end{tikzcd}
    \end{subfigure}%
    \begin{subfigure}{0.5\textwidth}
        \centering
        \begin{tikzcd}
            && {H^1} \\
            {H_1} && {H^1} \\
                  && {H^1} \\
            {H_2} && {H^2}
            \arrow["{B_2}", from=2-1, to=2-3]
            \arrow["{B_1}", curve={height=-12pt}, from=2-1, to=1-3]
            \arrow["{B_3}", from=4-1, to=3-3]
            \arrow["{B_4}", from=4-1, to=4-3]
        \end{tikzcd}
    \end{subfigure}
    \caption{A bipartite quiver $\mathcal{Q}$ and its constituent subspace quivers
    $\mathcal{Q}_1$ and $\mathcal{Q}_2$.}
    \label{F:splitting}
\end{figure}

\begin{proof}[Proof of sufficiency in \Cref{T:blqa}]
    By taking all but one subspace to be zero in condition \labelcref{blqadim}, we find that
    $\dim(V_i) \leq \sum_j \sum_{a \in \A_{ij}} p_j^{-1} \dim(B_a V_i)$ for all $V_i \subspeq H_i$ and
    each $i$. In particular, $\dim(H_i) \leq \sum_j \sum_{a \in \A_{ij}} p_j^{-1} \dim(H^j)$ for each $i$,
    so we also have $\dim(H_i) = \sum_j \sum_{a \in \A_{ij}} p_j^{-1} \dim(H^j)$ for each $i$ on account
    of condition \labelcref{blqascal}. Thus, the scaling condition \labelcref{scal} and the
    dimension condition \labelcref{dim} hold for each of the subspace quivers
    $\mathcal{Q}_i$ consisting of the source $H_i$, its incident arrows, and their targets
    \emph{regarded as separate vertices} (see \Cref{F:splitting} for an example), along with
    the corresponding weights $\vec{p}_i$.  It follows from \Cref{T:bl} that $\int_{H_i}
    \prod_j \prod_{a \in \A_{ij}} f_a \circ B_a \, dx \leq \BL(\mathcal{Q}_i,
    \vec{p}_i) \prod_j \prod_{a \in \A_{ij}}
    \norm{f_a}_{L^{p_j}(H^j)}$ with $\BL(\mathcal{Q}_i, \vec{p}_i) < \infty$ for each $i$;
    taking the product of these inequalities over $i$ yields the conclusion.
\end{proof}

Our proof that the conditions are necessary adapts an argument of Bennett, Carbery, Christ,
and Tao \cite{Bennett}.

\begin{proof}[Proof of necessity in \Cref{T:blqa}]
    Let $c$ be a positive constant to be determined later and let $r$ and $R$ be arbitrary
    positive constants satisfying $r \leq 1 \leq R$.  Given subspaces $V_i \subspeq H_i$,
    define
    \[
        S_i \defeq \set{(v_i, v_i') \in V_i \oplus V_i^\perp : 
        \norm{v_i} \leq c R \text{ and } \norm{v_i'} \leq c r} \subseteq H_i
    \]
    for each $i$ and
    \[
        S^a \defeq \set{(w_a, w_a') \in W^a \oplus (W^a)^\perp :
        \norm{w_a} \leq R \text{ and } \norm{w_a'} \leq r} \subseteq H^j
    \]
    for each $a \in \A_{ij}$ and each $i, j$, where $W^a \defeq B_a V_i \subspeq
    H^j$.

    If $(v_i, v_i') \in S_i$ for some $i$ and $a \in \A_{ij}$, then $B_a v_i \in W^a$ and
    $\norm{B_a v_i} \les \norm{v_i} \leq cR$, so $B_a v_i \in S^a$ provided that $c$ is chosen
    sufficiently small. 
    Similarly, if we write $B_a v_i'
    \eqdef (w_a, w_a') \in W^a \oplus (W^a)^\perp$, then 
    $\norm{w_a} + \norm{w_a'} \approx
    \norm{B_a v_i'} \les \norm{v_i'} \leq cr$, so $B_a v_i' \in S^a$ provided that
    $c$ is chosen sufficiently small because $r \leq R$.
    Hence $B_a S_i \subseteq S^a + S^a \subseteq 2S^a$.

    Now taking $f_a \defeq 1_{2S^a}$, we find that
    \begin{align*}
        \prod_{i} \int_{H_i} \prod_{j} \prod_{a \in \A_{ij}} f_a \circ
        B_a \, dx 
        &\geq
        \prod_{i} \int_{S_i} \prod_{j} \prod_{a \in \A_{ij}} f_a \circ
        B_a \, dx \\
        &=
        \prod_{i} \mes{S_i} \\
        &\approx
        R^{\sum_i \dim(V_i)} r^{\sum_i \codim(V_i)}.
    \end{align*}
    On the other hand,
    \begin{align*}
        \prod_i \prod_j \prod_{a \in \A_{ij}} \norm{f_a}_{L^{p_j}(H^j)}
        &\approx
        \prod_i \prod_j \prod_{a \in \A_{ij}} \mes{S^a}^{p_j^{-1}} \\
        &\approx
        R^{\sum_i \sum_j \sum_{a \in \A_{ij}} p_j^{-1} \dim(W^a)}
        r^{\sum_i \sum_j \sum_{a \in \A_{ij}} p_j^{-1} \codim(W^a)}.
    \end{align*}
    Sending $R \to \infty$, we deduce that 
    $\sum_i \dim(V_i) \leq \sum_i \sum_j \sum_{a \in \A_{ij}} p_j^{-1} \dim(W^a)$;
    sending $r \to 0$, we deduce that
    $\sum_i \codim(V_i) \geq \sum_i \sum_j \sum_{a \in \A_{ij}} p_j^{-1} \codim(W^a)$.
    Condition \labelcref{blqadim} is the first inequality, while condition
    \labelcref{blqascal} is the conjunction of the first inequality with $V_i \defeq H_i$
    and the second inequality with $V_i \defeq \set{0}$ (recalling that the $B_a$ are
    assumed to be surjective, so $B_a H_i = H^j$).
\end{proof}

\section{Examples of quiver Brascamp--Lieb inequalities} \label{S:examples}

Next, we return to the question of when an inequality of the form \labelcref{blq} holds,
which is answered by the work of Chindris and Derksen in the Gaussian case. Here we answer
this question in general by proving \Cref{T:equiv}, which implies that such an inequality is
equivalent to an inequality of the form \labelcref{blqa}.

\begin{proof}[Proof of \Cref{T:equiv}]
    Given inequality \labelcref{blqa} and a function $f_j$ for each $j$, we can take $f_a
    \defeq f_j$ for each $a \in \A_{ij}$ and each $i$ to obtain inequality
    \labelcref{blq} with constant $\BL(\mathcal{Q}, \vec{p})$, which shows that
    $\BLCD(\mathcal{Q}, \vec{p}) \leq \BL(\mathcal{Q}, \vec{p})$.
    
    Conversely, given inequality \labelcref{blq} and a function $f_a$ for each arrow $a$,
    let us assume that $\norm{f_a}_{L^{p_j}(H^j)} = 1$ for each $a$. We can then take $f_j
    \defeq \sum_i \sum_{a \in \A_{ij}} f_a$ for each $j$ to obtain
    \begin{align*}
        \prod_i \int_{H_i} \prod_j \prod_{a \in \A_{ij}} f_a \circ B_a \, dx 
        &\leq \BLCD(\mathcal{Q}, \vec{p})
        \prod_i \prod_j \prod_{a \in \A_{ij}} \Norm{\sum_{i'=1}^n \sum_{a' \in \A_{i'j}}
        f_a}_{L^{p_j}(H^j)} \\
        &\leq \BLCD(\mathcal{Q}, \vec{p})
        \prod_i \prod_j \prod_{a \in \A_{ij}} \alpha_j \\
        &= \BLCD(\mathcal{Q}, \vec{p}) \prod_j \alpha_j^{\alpha_j},
    \end{align*}
    which is inequality \labelcref{blqa} assuming that the $f_a$ are normalized.  By
    homogeneity, the inequality holds in general with the same constant, which means that
    $\BL(\mathcal{Q}, \vec{p}) \leq \prod_j \alpha_j^{\alpha_j} \cdot \BLCD(\mathcal{Q},
    \vec{p})$.
\end{proof}

Now let us see how \Cref{E:example} gives an instance of condition \labelcref{blqadim} being
strictly stronger than condition \labelcref{blgqdim}, thereby precluding a generalization of
Lieb's result (\Cref{T:lieb}) to the quiver setting.

\begin{proof}[Proof of \Cref{C:nolieb}]
    Consider the datum $(\mathcal{Q}, \vec{p})$ of \Cref{E:example}.  The scaling condition
    \labelcref{blgqscal} or \labelcref{blqascal} is $\dim(\R^3) = \frac{3}{2} \dim(\R^2)$,
    which is obviously satisfied. However, the dimension condition \labelcref{blgqdim} is
    $\dim(V_1) \leq \frac{3}{2} \dim(B_1 V_1 + B_2 V_1)$ for all $V_1 \subspeq \R^3$,
    whereas the dimension condition \labelcref{blqadim} is $\dim(V_1) \leq \frac{3}{4}
    \dim(B_1 V_1) + \frac{3}{4} \dim(B_2 V_1)$ for all $V_1 \subspeq \R^3$.  Evidently, the
    latter does not hold -- consider $V_1 = \spn \set{(1, 0, 0)}$ -- but in fact the former
    does. 

    To see this, first note that if $\dim(B_1 V_1 + B_2 V_1) = 0$, then $\dim(V_1) = 0$
    since $V_1 \subspeq \ker(B_1) \cap \ker(B_2) = \set{0}$, and if $\dim(B_1 V_1 + B_2 V_1)
    = 2$, the inequality is trivial. In the remaining case $\dim(B_1 V_1 + B_2 V_1) = 1$, we
    must have $B_1 V_1 \subspeq \spn \set{w}$ and $B_2 V_1 \subspeq \spn \set{w}$ for some
    $w = (w_1, w_2) \in \R^2$.  From this, we find that $V_1 \subspeq \spn \set{(w_1^2, w_1
    w_2, w_2^2)}$, so $\dim(V_1) \leq 1$ and the inequality is satisfied.

    As a result, $\BLCDG(\mathcal{Q}, \vec{p})$ is finite (by \Cref{T:blq}), yet
    $\BLCD(\mathcal{Q}, \vec{p})$ is not (by \Cref{T:blqa} or \Cref{C:equiv}).
\end{proof}

\begin{remark} \label{R:sum}
    For this datum, we can also observe directly that $\BLCD(\mathcal{Q}, \vec{p})$ is
    infinite: if $S \defeq ([0, N] \times [0, 1]) \cup ([0, 1] \times [0, N])$ and $f_1
    \defeq 1_S$, then the right-hand side of inequality \labelcref{blq} is
    $\norm{1_S}_{L^{4/3}(\R^2)}^2 = \mes{S}^{3/2} \approx N^{3/2}$, while the left-hand side
    is $\int_{\R^3} 1_S(x_1, x_2) 1_S(x_2, x_3) \, dx \geq \mes{[0, N] \times [0, 1] \times
    [0, N]} = N^2$.
\end{remark}

Even when $\BLCDG(\mathcal{Q}, \vec{p})$ and $\BLCD(\mathcal{Q}, \vec{p})$ are both finite,
they need not be equal -- in \Cref{E:R2-R1}, they are; in \Cref{E:R1-R1}, they are not (in
general).

\begin{example} \label{E:R2-R1}
    Consider the quiver $\mathcal{Q}$ in \Cref{F:example} with $H_1 \defeq \R^2$ and $H^1
    \defeq \R^1$, and let $\vec{p} = (p_1) \defeq (1)$. In addition, suppose that $B_1(x)
    \defeq b_1^\tp x$ and $B_2(x) \defeq b_2^\tp x$, where $b_1, b_2 \in \R^2$ are such that
    $B \defeq \begin{bmatrix} b_1 & b_2 \end{bmatrix}$ is invertible.

    According to the formula in \Cref{T:blq}, we have
    \begin{align*}
        \BLCDG(\mathcal{Q}, \vec{p})
        &=
        2 \cdot
        \sup_{a_1 > 0} \left[\frac{a_1^2}{\det(2b_1 a_1 b_1^\tp + 2b_2 a_1
        b_2^\tp)}\right]^{1/2} \\
        &=
        \frac{1}{\det(BB^\tp)^{1/2}} \\
        &=
        \frac{1}{\abs{\det(B)}}\,.
    \end{align*}
    On the other hand, for any measurable function $f_1 : \R^1 \to [0, \infty]$, we have
    \begin{align*}
        \int_{\R^2} f_1(b_1^\tp x) f_1(b_2^\tp x) \, dx
        &= \int_{\R^2} f_1(y_1) f_1(y_2) \, \abs{\det((B^\tp)^{-1})} \, dy \\
        &= \frac{1}{\abs{\det(B)}} \, \norm{f_1}_{L^1(\R^1)} \norm{f_1}_{L^1(\R^1)},
    \end{align*}
    which shows that
    \[
        \BLCD(\mathcal{Q}, \vec{p}) = \frac{1}{\abs{\det(B)}}
    \]
    as well.
\end{example}

\begin{example} \label{E:R1-R1}
    Consider the quiver $\mathcal{Q}$ in \Cref{F:example} with $H_1 \defeq \R^1$ and $H^1
    \defeq \R^1$, and let $\vec{p} = (p_1) \defeq (2)$. In addition, suppose that $B_1(x)
    \defeq b_1 x$ and $B_2(x) \defeq b_2 x$, where $b_1, b_2 \in \R^1 \setminus \set{0}$.

    According to the formula in \Cref{T:blq}, we have
    \begin{align*}
        \BLCDG(\mathcal{Q}, \vec{p})
        &=
        2^{1/2} \cdot
        \sup_{a_1 > 0} \left[\frac{a_1}{b_1 a_1 b_1 + b_2 a_1 b_2}\right]^{1/2} \\
        &=
        \left(\frac{2}{b_1^2 + b_2^2}\right)^{1/2}.
    \end{align*}
    On the other hand, for any measurable function $f_1 : \R^1 \to [0, \infty]$, we have
    \begin{align*}
        \int_{\R} f_1(b_1 x) f_1(b_2 x) \, dx
        &\leq \left[\int_{\R} f_1(b_1 x)^2 \, dx\right]^{1/2}
        \left[\int_{\R} f_1(b_2 x)^2 \, dx\right]^{1/2} \\
        &= \left[\int_{\R} f_1(y)^2 \, \abs{b_1^{-1}} \, dy\right]^{1/2}
        \left[\int_{\R} f_1(y)^2 \, \abs{b_2^{-1}} \, dy\right]^{1/2} \\
        &= \frac{1}{\abs{b_1 b_2}^{1/2}} \, \norm{f_1}_{L^2(\R^1)} \norm{f_1}_{L^2(\R^1)},
    \end{align*}
    which shows that
    \[
        \BLCD(\mathcal{Q}, \vec{p}) \leq \frac{1}{\abs{b_1 b_2}^{1/2}}\,.
    \]
    In fact, this is an equality. To see this, let $p > 1$ and
    \[
        f_1(x) \defeq
        \begin{cases}
            \abs{x}^{-p/2} & \text{if $\abs{x} \geq 1$}, \\
            0 & \text{otherwise}.
        \end{cases}
    \]
    Assuming without loss of generality that $\abs{b_1} \leq \abs{b_2}$, we compute that
    \begin{align*}
        \int_{\R} f_1(b_1 x) f_1(b_2 x) \, dx
        &= \int_{\abs{x} \geq 1/\abs{b_1}} \frac{1}{\abs{b_1 x}^{p/2}} \cdot
        \frac{1}{\abs{b_2 x}^{p/2}} \, dx \\
        &= \frac{1}{\abs{b_1 b_2}^{p/2}} \int_{\abs{x} \geq 1/\abs{b_1}}
        \left(\frac{1}{\abs{x}^{p/2}}\right)^2 \, dx \\
        &= \frac{1}{\abs{b_1 b_2}^{p/2}} \cdot \frac{1}{\abs{b_1}^{1-p}} \cdot \norm{f_1}_{L^2(\R^1)}^2,
    \end{align*}
    whence the claim follows by taking $p \to 1$.
    As a result, we see that both the Gaussian constant and the general constant are finite
    for this datum, but that they are equal if and only if $\abs{b_1} = \abs{b_2}$.
\end{example}

\section{Concluding remarks and questions}

Inspecting the proof of Lieb's theorem \cite{Lieb}, we find that it uses the
\emph{multilinearity} of inequality \labelcref{bl} in the $f_j$, which the quiver inequality
\labelcref{blq} does not possess. Thus, one might have expected the optimal constant for
Gaussians to differ from that for general functions in the latter inequality. Indeed, the
function $f_1$ in \Cref{R:sum} can be thought of as the sum of two rough approximations to
Gaussians, $g_1 \defeq 1_{[0, N] \times [0, 1]}$ and $h_1 \defeq 1_{[0, 1] \times [0, N]}$
(ignoring the overlap of their supports). The left-hand side of the inequality, roughly
$\int_{\R^3} (g_1 + h_1)(x_1, x_2) (g_1 + h_1)(x_2, x_3) \, dx$, is incomparably larger than
$\int_{\R^3} g_1(x_1, x_2) g_1(x_2, x_3) \, dx + \int_{\R^3} h_1(x_1, x_2) h_1(x_2, x_3) \,
dx$ for large $N$ because of the ``cross terms'' in the product, and consequently fails to
be bounded by the right-hand side.

Although the theory of quiver Brascamp--Lieb inequalities does not appear to be as rich as
that of ordinary Brascamp--Lieb inequalities, there are still some questions that could be
investigated. For instance, can the inequalities in \Cref{T:equiv} be strict? We know that
the second can be: consider the quiver $\mathcal{Q}$ in \Cref{F:example} once more, with
$H_1 \defeq \R^1$, $H^1 \defeq \R^1$, $B_1(x_1) \defeq x_1$, and $B_2(x_1) \defeq x_1$; and let
$\vec{p} = (p_1) \defeq (2)$. Then $\BLCD(\mathcal{Q}, \vec{p}) = \BL(\mathcal{Q}, \vec{p})
= 1$ by the Cauchy--Schwarz inequality and $\alpha_1 = 2$. We can also ask: if not
Gaussians, what are the maximizers of quiver Brascamp--Lieb inequalities? Are there any
interesting applications of such inequalities?

\emergencystretch=1em
\printbibliography

\end{document}